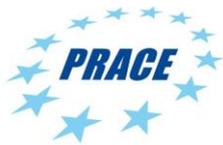

# A hybrid Hermitian general eigenvalue solver

Raffaele Solcà[*], Thomas C. Schulthess

*Institute forTheoretical Physics ETHZ, Swiss National Supercomputing Centre (ETHZ/CSCS)*

Azzam Haidar, Stanimire Tomov, Ichitaro Yamazaki, Jack Dongarra

*Institute Innovative Computing Laboratory, University of Tennessee*

**Abstract**

The adoption of hybrid GPU-CPU nodes in traditional supercomputing platforms opens acceleration opportunities for electronic structure calculations in materials science and chemistry applications, where medium sized Hermitian generalized eigenvalue problems must be solved many times. The small size of the problems limits the scalability on a distributed memory system, hence they can benefit from the massive computational performance concentrated on a single node, hybrid GPU-CPU system. However, new algorithms that efficiently exploit heterogeneity and massive parallelism of not just GPUs, but of multi/many-core CPUs as well are required. Addressing these demands, we implemented a novel Hermitian general eigensolver algorithm. This algorithm is based on a standard eigenvalue solver, and existing algorithms can be used. The resulting eigensolvers are state-of-the-art in HPC, significantly outperforming existing libraries. We analyze their performance impact on applications of interest, when different fractions of eigenvectors are needed by the host electronic structure code

## 1. Introduction

Scientific computing applications, ranging from computing frequencies that will propagate through a medium to the earthquake response of a bridge, or energy levels of electrons in nanostructure materials, require the solution of eigenvalue problems. There are many ways to formulate, mathematically, and solve these problems numerically [1]. In this work we are interested in dense eigensolvers, and in particular, generalized Hermitian-definite problems of the form $\mathbf{Ax} = \lambda\mathbf{Bx}$ where $\mathbf{A}$ is a Hermitian dense matrix and $\mathbf{B}$ is Hermitian positive definite. These solvers are needed for solving electronic structure problems in materials science and chemistry [2], [3], [4]. Particularly, problems with modest matrix dimensions of a few thousand to ten or twenty thousand seem to pose a challenge for most practical purposes. Typically these problems must be solved many times in the context of a parallel code. The known power limits in current processors, that would prevent the clock frequency to keep increasing with time, leaves us the challenge to solve the eigenvalue problem on nodes with a large number of threads that are executed on slower cores.

To implement a general eigensolver different subroutines are needed. Although these routines of interest are available in a number of standard libraries, current hardware changes have motivated their redesign (see Section 2) to achieve a more efficient use of the underlying hardware. In particular, hybrid architectures based on GPUs and multicore CPUs call for the development of new algorithms that can efficiently exploit the heterogeneity and the massive parallelism of not just the GPUs, but of the multi/many-core CPUs as well. Hybrid GPU-CPU nodes are already widely adopted in traditional supercomputing platforms such as the Cray XK6. This has opened many new acceleration opportunities in electronic structure calculations for materials science and chemistry applications, where generalized eigenvalue problems for medium sized Hermitian matrices must be solved many times, and they can benefit from the massive computational performance concentrated on a single node of a hybrid GPU-CPU system.

[*] Corresponding author. E-mail address: *rasolca@itp.phys.ethz.ch* .

## 2. Related Work

The main standard libraries that include eigensolver routines are LAPACK [5] and ScaLAPACK [6]. LAPACK is a linear algebra library for dense and banded matrices on shared-memory systems, while ScaLAPACK is its extension to distributed memory systems. LAPACK is based on calls to the BLAS (Basic Linear Algebra Subprograms) interface. Vendors usually provide well tuned LAPACK and ScaLAPACK versions, based on optimized BLAS implementations.

Work specifically on eigenvalue problems has been concentrated on accelerating separate components of the solvers and, in particular, the reduction of the matrix to tridiagonal form, which is the most time consuming task. The standard approach from LAPACK [7] is to use a "single phase" (also referred to as "one-stage") reduction. Alternatives are two- or more stage approaches, where the matrices are first reduced to band forms, and subsequently to their final tridiagonal form.

One of the first uses of a two-step reduction occurred in the context of out-of-core solvers for generalized symmetric eigenvalue problems [8]. With this approach, it was possible to recast the expensive memory-bound operations that occur during the panel factorization into a compute-bound procedure.

Consequently, a framework called Successive Band Reductions (SBR) was created [9], integrating some of the multi-stage work. The SBR toolbox applies two-sided orthogonal transformations based on Householder reflectors, and successively reduces the matrix bandwidth size until a suitable width is reached. The off-diagonal elements are then annihilated column-wise, which produces large fill-in blocks or bulges that need to be chased down towards the bottom right corner of the matrix. The bulge chasing procedure may result in substantial increase in the floating point operation count when compared with the standard single-phase approach from LAPACK. If eigenvectors are required in addition to eigenvalues, then the transformations may be efficiently accumulated using level 3 BLAS operations to generate these vectors. SBR heavily relies on multithreaded BLAS that are optimized to achieve satisfactory parallel performance. However, such parallelization paradigm incurs substantial overheads [10] as it fits the Bulk Synchronous Parallelism (BSP) model [11]. Communication bounds for such two-sided reductions have been established under the Communication Avoiding framework [12].

Tile algorithms have recently seen a rekindled interest when applied to the two-stage tridiagonal reductions [10]. Using high performance kernels combined with a data translation layer to execute on top of the tile data layout format, both implementations achieve a substantial improvement compared to the equivalent routines from the state-of-the-art numerical libraries. The off-diagonal elements are annihilated column-wise instead during the bulge chasing procedure, which engenders significant extra flops due to the size of the bulges introduced. An element-wise annihilation has then been implemented based on cache-aware kernels, in the context of the two-stage tridiagonal reduction [13]. Using the coalescing task technique [13], the performance achieved is by far greater than any other available implementation.

With the emergence of high bandwidth and high efficient devices such as GPUs, accelerating by orders of magnitude memory-bound and compute-bound operations became accessible. Tomov et al. [14] presented novel hybrid reduction algorithms for the two-sided factorizations, which take advantage of the high bandwidth of the GPU by offloading the expensive level 2 BLAS operations of the panel factorization to the GPU. The opposite was done by Bientinesi et al. [15] who accelerated the two-stage approach of the SBR toolbox by offloading the compute-bound kernels of the first stage to the GPU. The computation of the second stage (reduction to tridiagonal form) still remains on the host though.

Vomel et al. [16] extended the tridiagonalization approach in [14] to the symmetric standard eigenvalue problem.

A recent distributed-memory eigensolver library, developed for electronic structure codes, is ELPA [3]. ELPA is based on ScaLAPACK and does not support GPUs. It includes new one-stage and two-stages tridiagonalizations, the corresponding eigenvectors transformation, and a modified divide and conquer routine that can compute the entire eigenspace or just a part of it.

## 3. The general eigensolver

The steps to solve the general eigenproblem $\mathbf{Ax} = \lambda \mathbf{Bx}$ are described in Algorithm 1. The implementation of this algorithm combines the following routines:

- MAGMA Cholesky decomposition,
- MAGMA standard eigensolvers (Section 4),
- CUBLAS triangular solver,
- A hybrid implementation of the xHEGST routine that transforms the general problem into a standard problem, it is possible to build the general eigensolver.

The last routine has been implemented following the LAPACK algorithm, but moving all the level 3 BLAS operations are executed on the GPU.

---
1: Cholesky decomposition $\mathbf{B} = \mathbf{LL}^H$
2: transformation to standard eigenproblem $\mathbf{A}' = \mathbf{L}^{-1}\mathbf{AL}^{-H}$
3: standard eigensolver $\mathbf{A}'\mathbf{y} = \lambda \mathbf{y}$
4: eigenvectors backtransformation $\mathbf{x} = \mathbf{L}^{-H}\mathbf{y}$
---

Algorithm 1. Hermitian general eigenvalue solver

## 4. The standard eigensolver

The solution of the standard eigenproblem is described in Algorithm 2.

---
1: reduction to tridiagonal $\mathbf{T} = \mathbf{Q}_n^H \cdots \mathbf{Q}_2^H \mathbf{Q}_1^H \mathbf{A}' \mathbf{Q}_1 \mathbf{Q}_2 \cdots \mathbf{Q}_n$
2: tridiagonal eigensolver $\mathbf{Ty}' = \lambda \mathbf{y}'$
3: eigenvectors backtransformation $\mathbf{y} = \mathbf{Q}_1 \mathbf{Q}_2 \cdots \mathbf{Q}_n \mathbf{y}'$
---

Algorithm 2. Hermitian Standard eigensolver. *n* represents the number of stages in the tridiagonalisation.

### 4.1. One-stage standard eigensolver

The one-stage tridiagonalization is the approach used in the LAPACK eigensolvers. It was already implemented in MAGMA [16]. It has been improved including the hybrid algorithm of the tridiagonal eigensolver (Section 4.3) as well as including the possibility to compute only a fraction of the eigenspace. The tridiagonalisation's performances are limited, since 50% of the operations are level 2 BLAS operations.

### 4.2. Two-stage standard eigensolver

To solve the problem of the limited performance of the tridiagonalisation step, a two-stage approach has been used.

In the first step the matrix is reduced to a band matrix, while in the second step the band matrix is reduced to tridiagonal one using a bulge chasing technique. The algorithms are described in [13]. However, the reduction to band form is implemented using block techniques instead of tile techniques, since they are better for GPU operations. The algorithm is described in Figure 1.

For the second step we use a bulge chasing algorithm, very similar to [13]. However we differ by using a column-wise elimination, which allows us to have better locality for computing or applying the orthogonal matrix resulting from this phase. The drawback of this approach lies by the eigenvectors backtransformation, since the tridiagonalisation is performed in two steps, the backtransformation requires two steps too.

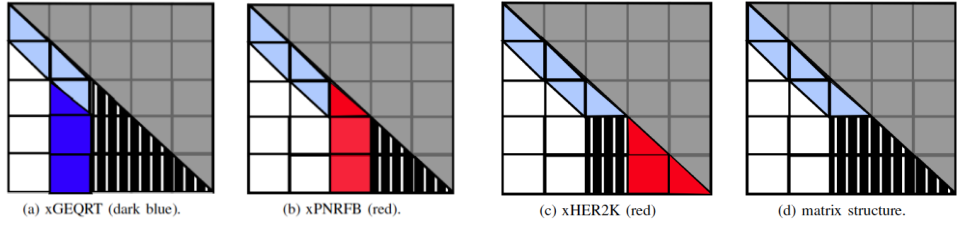

Figure 1. Reduction to band algorithm. The operation in (a) (CPUs) can be overlapped with (c) (GPU).

### 4.3. Tridiagonal eigensolver

The standard Divide and Conquer algorithm computes all the eigenvalues and eigenvectors of a real tridiagonal matrix. A detailed description of the algorithm can be found in [25]. This algorithm splits recursively the problem in two smaller problems. It then merges the solutions of the smaller problems as follows:

$$
\begin{aligned}
T &= \begin{pmatrix} T_1 & 0 \\ 0 & T_2 \end{pmatrix} + \rho \mathbf{v} \mathbf{v}^T \\
&= \begin{pmatrix} Z_1 & 0 \\ 0 & Z_2 \end{pmatrix} \left( \begin{pmatrix} \Lambda_1 & 0 \\ 0 & \Lambda_2 \end{pmatrix} + \rho \mathbf{u} \mathbf{u}^T \right) \begin{pmatrix} Z_1 & 0 \\ 0 & Z_2 \end{pmatrix}^T \\
(\mathbf{CPU}) &= \begin{pmatrix} Z_1 & 0 \\ 0 & Z_2 \end{pmatrix} Q \Lambda Q^T \begin{pmatrix} Z_1 & 0 \\ 0 & Z_2 \end{pmatrix} \\
(\mathbf{GPU}) &= Z \Lambda Z^T
\end{aligned}
$$

The multicore is used to compute the eigenvalues and eigenvectors of the rank-one modification, while the GPU updates the eigenvectors with a matrix-matrix multiplication.

It is easy to see that during the last merge one can compute only the needed eigenvectors, reducing the time to solution.

## 5. Results

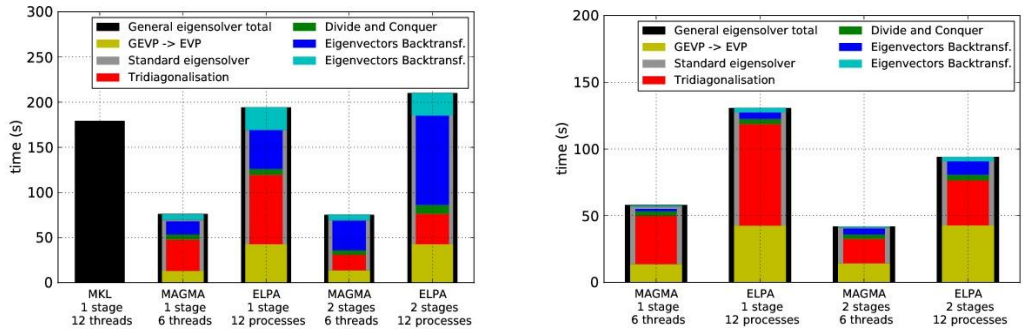

Figure 2. Time to solution of the double complex general eigensolver with matrix size 8000.
On the left the whole eigenspace is computed, while on the right only 10% of the eigenvectors.

Figure 2 shows the results of the hybrid algorithms compared with a shared memory library (MKL) and a distributed memory library (ELPA) with two different algorithms. To have a fair comparison on a 2 6-cores Intel Xeon 5650 and one Nvidia M2090 system we compare our hybrid routines, tested with six threads on one CPU socket and one GPU, against non-GPU routines tested using both of the CPU's sockets − 12 threads for shared-memory routines and 12 processes for distributed-memory processes.

The one-stage and two-stage approaches perform similarly for the complete eigenspace, while the two-stage approach shows better performance computing a fraction of the eigenspace. It every case the hybrid algorithms outperform the CPU- only algorithms by a factor more than 2.

## 6. MultiGPU support

We extended the dense symmetric generalized one-stage eigensolver of MAGMA to utilize multiple GPUs. The data on the GPUs is distributed in a 1D block cyclic way. Figure 3 shows the time spent in each step of the algorithm on up to eight Nvidia M2090 GPUs and two 6-cores Intel Xeon X5660.

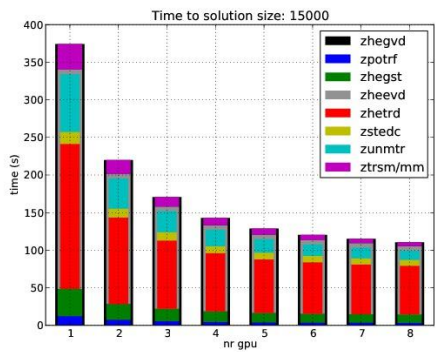

Figure 3. Time to solution for different number of GPUs of the double complex general eigensolver with matrix size 15000 and computing the whole eigenspace.

## 7. Conclusion

In this report we presented two different approaches to solve an eigenvalue problem benefiting from the computational power of a GPU. Both approaches show on a CPU socket+GPU a significant speedup with respect to shared-memory and distributed-memory approaches using two CPU sockets. The hybrid algorithms show more than a 2x speedup with respect to the CPUs-only libraries, and can be used to reduce the time needed for the electronic structure simulations, since the general eigensolver is the most time consuming part in these codes. Moreover, the extension of the routines to multiGPU allows a further reduction of the time to solution, as well as the handling of problems of larger size.

## Acknowledgements

This work was financially supported by the PRACE project funded in part by the EUs 7th Framework Programme (FP7/2007-2013) under grant agreement no. RI-211528 and FP7-261557. The work is partially achieved using the PRACE Research Infrastructure resources at Swiss National Supercomputing Centre (CSCS) in Switzerland.

The authors would also like to thank the National Science Foundation, the Department of Energy, NVIDIA, and the MathWorks for this research effort.